\documentclass[10pt]{amsart}            %%%  {article}

\usepackage{latexsym,amssymb}
\usepackage{amsmath,amsthm}      %%% for better alignments and theorems
\usepackage{amsopn}              %%% for `operator names' like `\det'
\usepackage{cite}                %%% for better handling of citations

\title{Differential forms and the Wodzicki residue.}

\author{William J. Ugalde}

\newcommand{\nnn}[1]{\eqref{#1}}

\newtheorem{theorem}{Theorem}[section]
\newtheorem{lemma}[theorem]{Lemma}
\newtheorem{proposition}[theorem]{Proposition}
\newtheorem{corollary}[theorem]{Corollary}
\newtheorem{definition}[theorem]{Definition}

\newcommand{\bbC}{\mathbb{C}}

\newcommand{\bbR}{\mathbb{R}}

\newcommand{\cA}{\mathcal{A}}
\newcommand{\cH}{\mathcal{H}}
\newcommand{\cL}{\mathcal{L}}
\newcommand{\cP}{\mathcal{P}}

\renewcommand{\a}{\alpha}
\renewcommand{\b}{\beta}
\renewcommand{\d}{\delta}
\newcommand{\ga}{\gamma}

\newcommand{\La}{\Lambda}
\newcommand{\Na}{\nabla}
\newcommand{\s}{\sigma}

\def\<#1,#2>{\langle{#1},{#2}\rangle} %% my version for inner product
\DeclareMathOperator{\tr}{trace}  %%% uses `amsopn' package

\def\set#1{\{\,#1\,\}}          %% set notation
\def\del{\partial}              %% partial derivatives symbol
\def\sideremark#1{\ifvmode\leavevmode\fi\vadjust{\vbox to0pt{\vss
\hbox to 0pt{\hskip\hsize\hskip1em
\vbox{\hsize2cm\tiny\raggedright\pretolerance10000
\noindent #1\hfill}\hss}\vbox to8pt{\vfil}\vss}}}

\def\lcovd#1#2{{#1}_{\,|\,#2}\,}        %% for f|_i
\def\ucovd#1#2{{#1}_{\,|\,}{}^{#2}\,}   %% for f|^i

\def\llcovd#1#2#3{{#1}_{\,|\,#2\,\,#3}\,}             %% for  f|_{ij}
\def\lucovd#1#2#3{{#1}_{\,|\,#2\,}{}^{#3}\,}          %% for  f|_i^j
\def\uucovd#1#2#3{{#1}_{\,|\,}{}^{#2\,#3}\,}          %% for  f|^{ij}

\def\ulucovd#1#2#3#4{{#1}_{\,|\,}{}^{#2\,}{}_{#3\,}{}^{#4}\,} %% for  f|^i_j^k
\def\llucovd#1#2#3#4{{#1}_{\,|\,#2\,#3\,}{}^{#4}\,}           %% for  f|_{ij}^k
\def\uuucovd#1#2#3#4{{#1}_{\,|\,}{}^{#2\,#3\,#4}\,}           %% for  f|^{ijk}
\def\lllcovd#1#2#3#4{{#1}_{\,|\,#2\,#3\,#4}\,}                %% for  f|_{ijk}

\renewcommand{\lcovd}[2]{{#1}_{\,;\,#2}}     %% f|_i
\renewcommand{\ucovd}[2]{{#1}_{\,;}{}^{#2}}  %% f|^i

\renewcommand{\llcovd}[3]{{#1}_{\,;\,#2#3}}         %% f|_{ij}
\renewcommand{\lucovd}[3]{{#1}_{\,;\,#2}{}^{#3}}    %% f|_i^j
\renewcommand{\uucovd}[3]{{#1}_{\,;}{}^{#2#3}}      %% f|^{ij}

\renewcommand{\ulucovd}[4]{{#1}_{\,;}{}^{#2}{}_{#3}{}^{#4}} %% f|^i_j^k
\renewcommand{\llucovd}[4]{{#1}_{\,;\,#2#3}{}^{#4}}         %% f|_{ij}^k
\renewcommand{\uuucovd}[4]{{#1}_{\,;}{}^{#2#3#4}}           %% f|^{ijk}
\renewcommand{\lllcovd}[4]{{#1}_{\,;\,#2#3#4}}              %% f|_{ijk}

\newcommand{\eps}{\varepsilon}      %% abbreviation
\renewcommand{\flat}{\mathrm{flat}} %% subscript for flat case
\DeclareMathOperator{\Hom}{Hom}   %% space of homomorphisms
\DeclareMathOperator{\Rc}{Rc}     %% Ricci tensor
\DeclareMathOperator{\Sc}{Sc}     %% scalar curvature
\DeclareMathOperator{\Wres}{Wres} %% Wodzicki residue
\DeclareMathOperator{\wres}{wres} %% Wodzicki residue density

\begin{document}

\begin{abstract}
For a pseudodifferential operator $S$ of order 0 acting on sections of a vector bundle $B$ on a compact manifold
$M$ without boundary, we associate a differential form of order dimension of $M$ acting on
$C^\infty(M)\times C^\infty(M)$.  This differential form $\Omega_{n,S}$ is given in terms of the
Wodzicki 1-density $\wres([S,f][S,h])$.  In the particular case of an even dimensional, compact,
conformal manifold without boundary,  we study this differential form for the case $(B,S)=(\cH,F)$,
that is, the Fredholm module associated by A. Connes \cite{Con1} to the manifold $M.$ We give its
explicit expression in the flat case and then we address the general case.  In \cite{Irazu} the author
presented the computations in the six-dimensional case of a whole family of differential forms related
to $\Omega_{n,S}.$
\par
\end{abstract}
%% 53A30 = Conformal differential geometry
%% 35S05 = General theory of PsDO
%% 58J40 = Pseudodifferential and Fourier integral operators on manifolds
%% 58B34 = Noncommutative geometry (à la Connes)
%% 46L87 = Noncommutative differential geometry

\maketitle  %%% Outputs the title, author(s) and date

\tableofcontents

\pagestyle{myheadings}
\markboth{William \ J. \ Ugalde}{Differential forms and the Wodzicki residue}

%%%%%%%%%%%%%%%%%%%%%%%%  INTRODUCTION %%%%%%%%%%%%%%%%%%%%%%%%%%%%%%%%%%%%%%%

\section{Introduction}

In \cite{Con1}, Connes associated to a even-dimensional compact
conformal manifold $M$ a canonical Fredholm module $(\cH,F)$ and use it
to define in the case $n=4$, an $n$-dimensional differential form $\Omega_n(f,h)$ on $M$
satisfying
\begin{equation}
\label{Wres=intOmega}
\Wres \bigl(f_0 [F,f][F,h]\bigr)=\int_M f_0 \Omega_n(f,h)
\end{equation}
for every $f_0,f,h$ in $C^\infty(M)$. This differential form is
determined by using the general formula for the total symbol of
the product of two pseudodifferential operators, it happens to be
\emph{conformally invariant} and \emph{uniquely} determined by the
relation \nnn{Wres=intOmega}. In the 4-dimensional case, this
differential form was found, by explicit computations, to be
\emph{symmetric} and such that the Paneitz operator $P$ relates to
it by $$ \int_M \Omega_4(f,h) =2 \int_M f P(h) d\,x. $$

Because of its potential relation with the study of conformally invariant differential operators
generalizing the Yamabe operator, like the \emph{GJMS operators} \cite{GJMS} for which the critical operator
$P_4=P/2$ in dimension $n=4$ was recovered from $\Omega_4$, the study of
such a differential form $\Omega_n$, is of considerable importance in
the case $n\geq 4$. In this paper we want to study this differential form in the following sense.

In Section~\ref{existsomega}, we associate to a pseudodifferential operator $S$ of order 0, acting on sections
of a vector bundle $B$ of rank $r$ on a compact manifold without boundary $M,$ an $n$-dimensional differential
form $\Omega_{n,S}$ acting on $C^\infty(M)\times C^\infty(M).$  $\Omega_{n,S}$ is obtained from the Wodzicki
1-density $\wres([S,f][S,h])$ by using the formula for the
total symbol of the product of two pseudodifferential operators \nnn{sympro}.
This differential form, will be uniquely determined
by the relation
$$
\Wres(f_0[S,f][S,h]) = \int_M f_0 \Omega_{n,S}(f,h)
$$
for every $f_0,f,h \in C^\infty(M)$.  Furthermore $\int_M f_0 \Omega_{n,S}(f,h)$ define
a Hochschild 2-cocycle on the algebra $C^\infty(M).$

In Section~\ref{omegaandfredmod} we will reduce ourselves to the
particular case $(B,S)=(\cH,F)$ of the Fredholm module associated
by Connes \cite{Con1} to an even-dimensional, compact, conformal
manifold without boundary. We will prove that this differential
form $\Omega_n = \Omega_{n,S}$ is indeed symmetric in $f$ and $h$
on any dimension. In Section~\ref{symbolofFflat} we compute the 
leading symbol of $F$, $\s_L^F$, and in Section~\ref{flatcase}, based on
$\tr\bigl(\s_L^F(\xi)\s_L^F(\eta)\bigr)$ for $\xi,\eta$ not zero
in $T^*_xM$ with the
operator under the trace sign acting on middle dimension forms we
will, in the flat case, express $\Omega_n$ as the integral of a
specific polynomial obtained from the Taylor expansion of the
function $\<\xi,\eta>^2|\xi|^{-2}|\eta|^{-2}$ around a diagonal
point. In Section~\ref{generalapproach}, we write $\Omega_n(f,h)$
as $(Q_n(df,dh) + Q_{R,n}(df,dh))d\,x$ for a universal bilinear
form which coincides with $\Omega_n$ in the flat case  (see
\cite{Con1} for the case $n=4$), and then we address the
computations for the general case.  At the end, we present an
appendix on the dependence of the symbol of $F$ on the partial
derivatives of the metric.

In \cite{Irazu} the author presented a first glance at the
computations of this differential form in the 6-dimensional case.
The results on this paper are the rigorous proof of the statements
needed for that computations, and any other done in general
dimensions, to be true.  My special thanks to T. Branson for
unlimited support,  and to J. C. V\'arilly for helpful
comments on an earlier version.

%%%%%%%%%%%%%%%%%%%%%%%%%%%%  SETTING AND NOTATIONS   %%%%%%%%%%%%%%%%%%%%%%%%%%%%%%%

\medskip
\noindent{\bf Setting and Notations}

Our general setting comes from \cite{Con1} where, using the quantized
calculus, Connes gives an analogue in dimension 4 of the 2-dimensional
Polyakov action
$$
I(X) = \frac{1}{2\pi} \int_\Sigma \eta_{ij} \,dX^i \wedge \star dX^j
$$
for a Riemann surface $\Sigma$ and a map $X$ from $\Sigma$ to $\bbR^d$.
In the present work, unless otherwise stated, $M$ represents an oriented,
compact, manifold without boundary, endowed with a conformal structure,
and of even dimension $n=2l$ (which is going to be fixed from now on).
$\cH$ is the Hilbert space of square integrable forms of middle
dimension, $\cH = L^2(M,\Lambda^{l}_{\bbC} T^*M)$, on which functions
on $M$ act as multiplication operators.
$F$ is the pseudodifferential operator of order 0 acting in $\cH,$
obtained from the orthogonal projection $P$ on the image of $d,$ by the relation $F=2P-1.$
From the Hodge decomposition theorem \cite{Warner}
it is easy to see that $F$ preserves the finite
dimensional space of harmonic forms $H^l$, and that $F$ restricted to the $\cH\ominus H^l$
is given by
\begin{equation}
\label{defofF}
F=\frac{d \d-\d d}{d \d+\d d},
\end{equation}
in terms of a Riemannian metric compatible with the conformal structure of $M.$
Both $\cH$ and $F$ are independent of the metric in the conformal class \cite[Section~IV.4.$\ga$]{Koran}.

The Riemann curvature tensor will be represented with the letter $R$, the Ricci tensor will be
represented by $\Rc_{ij} = R^{k}{}_{ikj}$, and the scalar curvature by $\Sc = \Rc^{i}{}_{i}$.
If needed, we will ``raise'' and ``lower''
indices using the metric without explicit mention, for example, $g_{mi}R^i{}_{jkl}=R_{mjkl}$.

When working with the total symbol of a pseudodifferential operator
$P,$ we will denote its leading symbol by $\s_L^P,$ or $\s_L(P)$ in case $P$ has a long expression,
and if the operator $P$ is of order $k$ then its total symbol (in some given local coordinates)
will be represented as
$$
\s(P) = \s_k^P + \s_{k-1}^P + \s_{k-2}^P+\cdots,
$$
where $\s_L^P=\s_k^P.$  It is important to note that the different $\s_j^P$ for $j<k$ are defined
only in local charts and are not diffeomorphism invariant \cite{Ku}.  However, Wodzicki~\cite{Wo} has shown
that the term $\s_{-n}^P$ enjoys a very special significance.   For a pseudodifferential operator $P,$ acting
on sections of a vector bundle $B$ over a manifold $M,$ there is a $1$-density on $M$ expressed in local
coordinates by
\begin{equation}
\label{wresP}
\wres(P) = \int_{|\xi|=1} \left\{\tr(\s_{-n}^P(x,\xi))\,d^{n-1}\xi\right\} \,d\,x.
\end{equation}
This \emph{Wodzicki residue density} is independent of the local representation.  Here we are using the same
notations as in \cite{Polaris}, where an elementary proof of this matter can be found.

In the case of a Riemannian manifold $(M,g)$ and an orientation $d\,x$, the metric $g$ induces a
metric
$g^1$ on the cotangent bundle $T^* M$, and metrics $g^k$ on the exterior bundles $\La^k M$ in the following
way \cite{Bra3}.  If $X$ is a vector field and $\omega_X$ denotes the 1-form determined by
$\omega_X(Y)=g(X,Y)$
then $g^1(\omega_X,\omega_Y):=g(X,Y)$.  If $\omega^i,\eta^j$ are 1-forms then,
$$
g^k(\omega^1\wedge\cdots\wedge\omega^k,\eta^1\wedge\cdots\wedge\eta^k):=\det(g^1(\omega^i,\eta^j)).
$$
Saying that $d\,x$ is normalized means it has been chosen so that $g^n(d\,x,d\,x)=1$.

%%%###########%%% END of NOTATIONS %%%###########%%%

%%%%%%%%%%%%%%%%%%%%%%%%section: Properties of $\Omega_{n,S}$ %%%%%%%%%%

\section{Properties of $\Omega_{n,S}$}

For a pseudodifferential operator $S$ of order 0 acting on sections of a vector bundle $B$ over a
compact manifold without boundary $M$,
we consider $P$ to be the pseudodifferential operator of order $-2$ given by the product
$P= f_0[S,f][S,h]$ with each $f_0,f,h \in C^\infty(M)$.
$P$ is acting on the same vector bundle as $S$, where smooth functions on $M$ act as multiplication operators.
The total symbol of $P$, up to order $-n$,
is represented as a sum of $r \times r$ matrices of the form
$$
\s_{-2}^P+\s_{-3}^P+\cdots+\s_{-n}^P
$$
where $r$ is the rank of the vector bundle $B.$
We aim to study
$$
\Wres(P) = \int_M \left\{\int_{|\xi|=1}\tr(\s^P_{-n}(x,\xi))\,d^{n-1}\xi\right\}\,d\,x
$$
where $\s_{-n}(x,\xi)$ is the component of order $-n$ in the total symbol of $P,$
$|\xi|=1$ means the Euclidean norm of the coordinate vector
$(\xi_1,\cdots,\xi_n)$ in $\bbR^n,$ and $d^{n-1}\xi$ is the normalized volume on
$\set{|\xi|=1}.$
$\Wres(P)$ is independent of the choice
of the local coordinates on $M$, the local basis of $B$, and defines a trace (see \cite{Wo}).

%%%%%%%%%%%%%%%%%%%%%%%%section: Existence of $\Omega_{n,S}$ %%%%%%%%%%

\subsection{Existence and uniqueness of $\Omega_{n,S}$}
\label{existsomega}

In general, the total symbol of the product of two pseudodifferential
operators $P_1$ and $P_2$ is given by
\begin{equation}
\label{sympro}
\sigma(P_1 P_2) = \sum \frac{1}{\a!}
\del^\a_\xi(\sigma(P_1)) \,D^\a_x(\sigma(P_2))
\end{equation}
where $\a=(\a_1,\dots,\a_n)$ is a multi-index,
$\a! = \a_1!\cdots \a_n!$ and $D^\a_x = (-i)^{|\a|} \del^\a_x$.

\begin{lemma}
Let $T$ be a pseudodifferential operator and for $f$ in $C^\infty(M)$ let $f$
represent the operator ``multiplication by $f$'' then $\s(fT) = f\s(T)$
and in particular $\wres(fT) = f \wres(T)$, implying
\begin{equation}
\label{wresf0P}
\wres(f_0[S,f][S,h]) =
f_0 \wres([S,f][S,h]).
\end{equation}
\end{lemma}

\begin{proof}
By \nnn{sympro} $\s(f T)=\s(f)\s(T)$
since $\s(f)$ does not depend on $\xi$, indeed $\s(f) = f I$ with $I$ the identity
operator on the sections the operator $T$ is acting on,
and the first assertion follows.
In particular note that $\s_{-n}(f T) = f \s_{-n}(T)$.

By \nnn{wresP} we have
\begin{align*}
\wres(fT)
&= \left\{\int_{|\xi|=1} \tr(\s_{-n}(f T) \,d^{n-1}\xi)\right\} \, d \,x \\
&= \left\{\int_{|\xi|=1} f \tr(\s_{-n}(T) \,d^{n-1}\xi)\right\} \, d \,x = f \wres(T),
\end{align*}
and the lemma follows by taking $T=[S,f][S,h]$.
\end{proof}

\begin{lemma}
$[S,f]$ is a pseudodifferential operator of order $-1$ with total symbol
$\s([S,f]) = \sum_{k \geq 1} \s_{-k}([S,f])$ where
\begin{equation}
\label{sym-nSf}
\s_{-k}([S,f]) =\sum_{|\b|=1}^k \frac{1}{\b!} D^\b_x(f)
     \del_\xi^\b (\s^S_{-(k-|\b|)}).
\end{equation}
\end{lemma}

\begin{proof}
Since $S$ and $f$ are of order zero, the commutator $[S,f]$ is of order minus one.
The expansion of
the symbol $\s([S,f])$ is given by
\begin{align*}
\s([S,f])
&= \sum_{|\b|\geq 1} \frac{1}{\b!} \bigl(\del_\xi^\b(\s^S)\,D^\b_x(f)\bigr) \\
&= \sum_{|\b|=1} \frac{1}{\b!} D^\b_x(f)
   \bigl(\del_\xi^\b (\s^S_0) + \del_\xi^\b (\s^S_{-1}) +
   \del_\xi^\b (\s^S_{-2}) + \cdots \bigr) \\
&\qquad \qquad + \sum_{|\b|=2} \frac{1}{\b!} D^\b_x(f)
     \bigl(\del_\xi^\b (\s^S_0) + \del_\xi^\b (\s^S_{-1})  + \cdots \bigr) \\
&\qquad \qquad \qquad \qquad+ \sum_{|\b|=3} \frac{1}{\b!} D^\b_x(f)
     \bigl(\del_\xi^\b (\s^S_0) + \cdots \bigr)
     + \cdots \\
&= \sum_{k\geq 1} \biggl(\sum_{|\b|=1}^k \frac{1}{\b!} D^\b_x(f)
     \del_\xi^\b (\s^S_{-(k-|\b|)}) \biggr),
\end{align*}
and the result follows.
\end{proof}

\begin{lemma}
\label{lemmas_nSS}
With the sum taken over $|\a'|+|\a''|+|\b|+|\d|+j+k = n$, $|\b|\geq 1$,
and $|\d|\geq 1$,
\begin{align}
&\s_{-n}([S,f][S,h])
\nonumber \\
&=\sum \frac{1}{\a'!\a''!\b!\d!}
   D^\b_x(f) D_x^{\a''+\d}(h) \del_\xi^{\a'+\a''+\b}(\s^S_{-j})
\del_\xi^\d(D_x^{\a'}(\s^S_{-k})).
\label{sigma-nSS}
\end{align}
\end{lemma}

\begin{proof}
By \nnn{sympro} we have
$$
\s_{-n}([S,f][S,h])
= \sum_{r+s+|\a|=n}\frac{1}{\a!}
   \del^\a_\xi(\s_{-r}([S,f])) D^\a_x(\s_{-s}([S,h])).
$$
Using \nnn{sym-nSf} and because $\del^\a_\xi(D^\b_x f) = 0$ for any $|\a|\geq 1$,
we deduce
\begin{align*}
&\s_{-n}([S,f][S,h]) \\
&=
 \sum_{{r+s+|\a|=n}\atop{r \geq1,s \geq1}} \frac{1}{\a!}
   \sum_{{j+|\b|=r}\atop{|\b|\geq1}} \frac{1}{\b!}
     \del^\a_\xi(D^\b_x(f) \del^\b_\xi(\s^S_{-j}))
  \sum_{{k+|\d|=s}\atop{|\d|\geq1}} \frac{1}{\d!}
     D^\a_x(\del^\d_\xi(\s^S_{-k}) D^\d_x(h)) \\
&
= \sum_{{r+s+|\a|=n}\atop{r \geq1,s \geq1}} \frac{1}{\a!}
   \sum_{{j+|\b|=r}\atop{|\b|\geq1}} \frac{1}{\b!}
     D^\b_x(f) \del^{\a+\b}_\xi(\s^S_{-j})
\times \\
&\qquad \times \,
   \sum_{{k+|\d|=s}\atop{|\d|\geq1|}} \frac{1}{\d !}
   \sum_{{0\leq\a' \leq \a}\atop{\a'+\a''=\a}} \frac{\a!}{\a'!\a''!}
     \del^\d_\xi(D^{\a'}_x(\s^S_{-k})) D^{\a''+\d}_x(h),
\end{align*}
where we use, in the last equality, the Leibniz rule:
$$
D^\a_x(fh) = \sum_{{0 \leq \a' \leq \a}\atop{\a'+\a''=\a}}
\frac{\a!}{\a'!\a''!} D^{\a'}_x(f) D^{\a''}_x(h),
$$
with $\a'$ a multi-index of the form $\a'=(\a'_1,\dots,\a'_n)$, and
$\a'\leq\a$ meaning $\a'_i \leq \a_i$ for every $i=1,\dots,n$, the sum
taken over $|\a|+|\b|+|\d|+j+k = n$, $|\b|\geq 1$, and $|\d|\geq 1$.
\end{proof}

By \nnn{wresf0P},
\begin{align}
\wres(f_0[S,f][S,h]))
&= f_0 \biggl\{ \int_{|\xi|=1}
   \tr\biggl\{ \sum \frac{1}{\a'!\a''!\b!\d!}
   D^\b_x(f) D^{\a''+\d}_x(h)
\times \nonumber \\
&\qquad \times
   \del^{\a'+\a''+\b}_\xi(\s^S_{-j}) \del^{\d}_\xi(D^{\a'}_x(\s^S_{-k}))
\biggr\}\,d^{n-1}\xi \biggr\}\,d\,x
\nonumber \\
&= f_0 \omega_{n,S}(f,h)\, dx = f_0 \Omega_{n,S}(f,h).
\label{wresf0f1f2}
\end{align}
which explicitly gives an expression for $\Omega_{n,S}(f,h)$.

\begin{definition}
For every $f,h \in C^\infty(M)$ we define by \nnn{wresf0f1f2}
$$
\Omega_{n,S}(f,h) =\omega_{n,S}(f,h)\,d\,x := \wres([S,f][S,h]).
$$
\end{definition}

Because of the previous construction, we have the $\Omega_{n,s}$ is uniquely determined by its relation with the
Wodzicki residue of the operator $f_0[S,f][S,h]$ as stated in the following theorem.

\begin{theorem}
\label{theoremonOmegan1}
There is a unique, $n$-diffe\-ren\-tial
form $\Omega_{n,S}$ such that
$$
\Wres(f_0[S,f][S,h]) = \int_M f_0 \Omega_{n,S}(f,h)
$$
for all $f_0,$ $f,$ $h$ in $C^\infty(M)$.
\end{theorem}

By construction, $\Wres(f_0[S,f_1][S,f_2])$ is a Hochschild 
2-cocycle over the algebra of smooth functions on
a compact manifold, $\cA = C^\infty(M)$. Next we want to state such a property.
A Hochschild $n$-cochain $\varphi$ on $\cA$ is
an $(n+1)$-linear functional on $\cA$. The coboundary
operator, denoted by $b$, is given by:
\begin{align*}
&(b \varphi) (f_0,\dots,f_{n+1}) := \\
& \sum_{j=0}^n (-1)^j \varphi(f_0,\dots,f_jf_{j+1},\dots,f_{n+1}) -
 (-1)^{n+1} \varphi(f_{n+1}f_0,\dots,f_n).
\end{align*}
The cohomology of this complex is the {\bf Hochschild cohomology} of
$\cA$. In particular, a Hochschild 0-cocycle $\varphi$ on $\cA$ is a
trace since $\varphi \in \cA^* =\Hom(\cA,\bbC)$, and
$0 = (b \varphi)(f_0,f_1) = \varphi(f_0f_1) - \varphi(f_1 f_0).$
A Hochschild 2-cocycle $\varphi$ must satisfy
\begin{align*}
0 &= (b\varphi)(f_0,f_1,f_2,f_3) \\
  &= \varphi(f_0 f_1,f_2,f_3) - \varphi(f_0,f_1 f_2,f_3) +
     \varphi(f_0,f_1,f_2 f_3) - \varphi(f_3 f_0,f_1,f_2)
\end{align*}
for every $f_i \in C^\infty(M).$

\begin{lemma}
\label{WresHochschild2cocycle}
$\Wres(f_0[S,f_1][S,f_2])$ is a Hochschild 2-cocycle over the 
algebra of smooth functions on $M.$
\end{lemma}

\begin{proof}
Linearity is evident, and from the relation $[S,fh] = [S,f]h + f[S,h]$
we have, for $\varphi(f_0,f_1,f_2) = \Wres(f_0[S,f_1][S,f_2])$:
\begin{align*}
(b\varphi)(f_0,f_1,f_2,f_3)
&= \Wres(f_0 f_1[S,f_2][S,f_3]) - \Wres(f_0[S,f_1 f_2][S,f_3]) \\
&\qquad + \Wres(f_0[S,f_1][S,f_2 f_3]) - \Wres(f_3 f_0[S,f_1][S,f_2]) \\
&= \Wres(f_0 f_1[S,f_2][S,f_3]) \\
&\qquad - \Wres(f_0[S,f_1]f_2[S,f_3]) - \Wres(f_0 f_1[S,f_2][S,f_3]) \\
&\qquad + \Wres(f_0[S,f_1][S,f_2]f_3) + \Wres(f_0[S,f_1]f_2[S,f_3]) \\
&\qquad - \Wres(f_3 f_0[S,f_1][S,f_2]) \\
&= \Wres(f_0[S,f_1][S,f_2]f_3) - \Wres(f_3 f_0[S,f_1][S,f_2]) = 0
\end{align*}
because of the trace property of $\Wres$.
\end{proof}

So far, by taking $f_0 = 1$, by uniqueness, and by the trace
property of $\Wres$, we conclude $$ \int_M \Omega_{n,S}(f,h) =
\Wres([S,f][S,h]) =  \Wres([S,h][S,f]) = \int_M \Omega_{n,S}(h,f).
$$ From \nnn{wresf0P}, $\Omega_{n,S}(f,h) = \wres([S,f][S,h])$ but
in general $\wres$ is not a trace, hence asserting that
$\Omega_{n,S}(f,h)$ is symmetric is asserting that $$
\wres([S,f][S,h])=\wres([S,h][S,f]). $$ To conclude the symmetry
of $\Omega_{n,S}(f,h)$ on $f$ and $h$, it is necessary to request
more properties on the operator $S$.  As we will see, it is enough
to have the property $S^2f=fS^2,$ for every $f\in C^\infty(M).$

\begin{theorem}
\label{symmetryofOmega}
If $S^2 f = f S^2$ for every $f \in
C^\infty(M)$ then the differential form $\Omega_{n,S}$ associated
by Theorem~\ref{theoremonOmegan1} to $(B,S)$ is symmetric in $f$
and $h$.
\end{theorem}

\begin{proof}
We are going to exploit the linearity and the trace property of
the Wodzicki's residue \cite{Wo}.  Note that
\begin{align*}
&\Wres(f_0[S,f][S,h] - f_0[S,h][S,f])
\\
&=\Wres(f_0 S f S h -f_0S
f h S - f_0 f S S h + f_0 f S h S
\\
&\quad\qquad - f_0 S h S f +
f_0 S h f S + f_0 h S S f - f_0 h S f S)
\\
&=\Wres(f_0 S f S h +
f_0 f S h S - f_0 S h S f - f_0 h S f S)
\end{align*}
using the fact $S^2$ commutes with any element of the algebra
$C^\infty(M)$ and the commutativity of the algebra $C^\infty(M)$.
By the trace property of the Wodzicki's residue, it follows that
$\Wres(f_0 F f F h)=\Wres(h f_0 F f F)$ and $\Wres(f_0 F h F
f)=\Wres(f f_0 F h F)$. Hence, by using once more the
commutativity of $C^\infty(M)$ we conclude $$ \int_M f_0
\Omega_{n,S}(f,h) = \int_M f_0 \Omega_{n,S}(h,f) $$ for every
$f_0,f,h \in C^\infty(M)$. Therefore $\Omega_{n,S}(f,h) =
\Omega_{n,S}(h,f)$ because of the arbitrariness of $f_0$.
\end{proof}

%%%%%%%%%%%%%%%%%%%%  subsection:  $\Omega_{n,S}$ and the Fredholm module on conformal manifolds  %%%%%%%%%%%%%

\subsection{$\Omega_{n,S}$ and the Fredholm module on conformal manifolds}
\label{omegaandfredmod}

For the rest of this paper we restrict ourselves to an even
dimensional, compact, oriented, conformal manifold without
boundary $M$, and $(B,S)$ given by the canonical Fredholm module
$(\cH,F)$ associated to $M$ by A. Connes \cite{Con1}. In this
particular case, the pseudodifferential operator of order 0 is
given by $F = (d \d - \d d)(d \d + \d d)^{-1}$ acting on $\cH =
L^2(M,\La^{l}_\bbC T^*M)\ominus H^l,$ with $H^l$ the finite
dimensional space of middle dimension harmonic forms.   By
definition $F$ is selfadjoint and such that $F^2=1.$ We relax the
notation by denoting $\Omega_n = \Omega_{n,F}$ in this particular
situation.

\begin{theorem}
\label{theoremonOmegan}
In the particular case in which $M$ is a even dimensional compact conformal manifold without boundary and
$(\cH,F)$ is the Fredholm module associated to $M$ by A. Connes \cite{Con1},
there is a unique, symmetric, and conformally invariant $n$-diffe\-ren\-tial form $\Omega_n=\Omega_{n,F}$
such that
$$
\Wres(f_0[F,f_1][F,f_2]) = \int_M f_0 \Omega_n(f_1,f_2)
$$
for all $f_i \in C^\infty(M).$  Furthermore, $\int_M f_0 \Omega_n(f_1,f_2)$ defines a Hochschild 2-cocycle
on the algebra of smooth functions on $M.$
\end{theorem}

\begin{proof}
Uniqueness follows from \nnn{wresf0f1f2}.  Symmetry follows from
Theorem~\nnn{symmetryofOmega} and its conformal invariance follows
from its construction. Indeed, as pointed out in
\cite[Section~2.1.2]{A-M}, $\Wres(P)$ does not depend on the
choice of the metric in the conformal class defining the cosphere
bundle.  The only possible metric dependence is given by the
operator $P$.  In our particular case $P=[F,f][F,h]$ is, as well
as $F$, independent of the metric in the conformal class. Indeed,
the only metric dependence of $(d \d -\d d)(d \d + \d d)^{-1}$ is
the one given by the Hodge star operator $\star$ used in $\d =
-\star d \star$, which is invariant under conformal changes of the
metric when restricted to middle dimension forms (see
Lemma~\nnn{hodgestar}).  The last assertion follows from
Lemma~\nnn{WresHochschild2cocycle}.
\end{proof}

\begin{lemma}
\label{hodgestar}
In the particular case of forms of middle dimension, the restriction
$\star|_{\La^{l}(T^*M)}$ is conformally invariant.
\end{lemma}

\begin{proof}
If $\hat g = e^{2 \eta}g$ is a conformally equivalent metric to
the metric $g$, then on $k$-forms we have ${\hat g}^k = e^{-2 k
\eta} g^k$. In particular, with $\hat{d \,x}$ and $\hat{\star}$
the volume form and the Hodge star operator associated to the
metric $\hat{g}$ we have $$ 1 = {\hat{g}}^n(\hat{d
\,x},\hat{d\,x}) = e^{-2 n \eta} \,g^n(\hat{d \,x},\hat{d \,x}).
$$ Now for two $(n-k)$-forms $\xi$ and $\xi'$ with $\xi \wedge
\star\xi' = g^{n-k}(\xi,\xi') \,d \,x$, we get $$ \xi \wedge
\hat{\star} \xi' = {\hat g}^{n-k}(\xi,\xi') \hat{d \,x}
  = e^{-2(n-k)\eta} g^{n-k}(\xi,\xi') e^{n \eta} d \,x
  = e^{(2 k-n)\eta} \xi \wedge \star\xi'.
$$
and hence $\hat{\star} = e^{(2 k-n)\eta}\star$. This shows that for
$k = l$, the scaling has had no effect.
\end{proof}

%%%%%%%%%%%%%%%%%%%%%%%%%%%%%%%%%%%%%%%%%%%%%%%%%%%%%%%%%%%%%%%%%%%%

%%%%%%%%%%%%%%% section : The leading symbol of $F$ %%%%%%%%%%%%%%%

\section{The leading symbol of $F$}
\label{symbolofFflat}

In this section we study the symbol of $F$ in the flat case, and in
general, its leading symbol.
In \cite[Lemma~1.5.3]{Gilkey}, it is proved that the leading symbols of
$d$ and $\d$ are given by $\s_L(d) = i \eps(\xi)$, $\s_L(\d) = -i \iota(\xi)$, where
$\xi \in T^*M$, and $\eps(\xi)$ and $\iota(\xi)$ are the exterior and interior multiplications
respectively.  As a result of the identity
$(i \eps(\xi) - i \iota(\xi))^2 = |\xi|^2 I$
and since $\Delta = d \d + \d d = (d+\d)^2$, in the same reference it
is also proved that the leading symbol of $\Delta$ is given by
$\s_L(\Delta) = |\xi|^2 I.$

Because of the rule for the leading symbol of the product of two
pseudodifferential operators, $\s_L(P_1 P_2) = \s_L(P_1) \s_L(P_2)$,
we conclude that $\s_L(\Delta^{-1}) = |\xi|^{-2} I$.

For the particular case of a flat metric, the total
symbols of both $d \d$ and $\d d$ coincide with their leading symbols,
$\s(d \d) = \s_L(d \d)$ and $\s(\d d) = \s_L(\d d)$, which at the same
time, do not depend on the variable $x$. As a consequence, the total
symbols of both $\Delta = d \d + \d d$ and $d \d - \d d$ coincide with
their leading symbols $\s(\Delta) = \s_L(\Delta)$ and
$\s(d \d-\d d) = \s_L(d \d-\d d)$ in the flat case.
It follows also that
---see Lemma~\ref{lemasymdelta-1}--- the total symbol of $\Delta^{-1}$
coincides with its leading symbol $\s(\Delta^{-1}) = \s_L(\Delta^{-1})=|\xi|^{-2}I$
and furthermore, it does not depend on $x$, hence
\begin{align*}
\s(F)(x,\xi)
&= \s(d \d-\d d)(x,\xi) \,\s(\Delta^{-1})(x,\xi) \\
&= |\xi|^{-2} (\eps(\xi)\iota(\xi) - \iota(\xi)\eps(\xi)).
\end{align*}
Summarizing, we have:

\begin{proposition}
\label{leadingsymbolofF}
The leading symbol of $F$ is given by
$$
\s_L(F)(x,\xi) = \s_L(F)(\xi) =
|\xi|^{-2}(\eps_\xi \iota_\xi - \iota_\xi \eps_\xi)
$$
for all $(x,\xi) \in T^* M$, $\xi \not= 0$. In the particular case of
a flat metric, we also have $\s_{-k}(F) = 0$ for all $k \geq 1$.
\end{proposition}

%%%%%%%%%%%%%%%%%%%%%%%%%%%%%%%%%%%%%%%%%%%%%%%%%%%%%%%%%%%%%%%%%%%%%%

%%%%% section: $\Omega_n$ in the flat case %%%%%

\section{$\Omega_n$ in the flat case}
\label{flatcase}

Let $\Omega_{n \,\flat}(f,h)$ be the $n$-dimensional form on a flat
manifold $M$ uniquely determined by Theorem~\ref{theoremonOmegan1}.
Since we are in the flat case we have
$\s_{-k}^F(x,\xi) = 0$,
for all $k > 0$.  Because of that and \nnn{sym-nSf},
$$
\s_{-r}([F,f]) = \sum_{|\b|=r} \frac{1}{\b!} \del_\xi^\b(\s_L^F) D_x^\b f.
$$
Using this information we deduce from \nnn{sigma-nSS} with $\a' = 0$ and
$\a'' = \a$,
$$
\s_{-n}([F,f][F,h]) =
\sum \frac{1}{\a!\b!\d!}
   (D^\b_x f)(D^{\a+\d}_x h)(\del^{\a+\b}_\xi(\s_L^F))(\del^\d_\xi(\s_L^F)),
$$
with the sum taken over $|\a|+|\b|+|\d| = n$, $1 \leq |\b|$, $1\leq |\d|$.
By \nnn{wresf0f1f2} and Theorem~\ref{theoremonOmegan1}, we have:

\begin{lemma}
\begin{align*}
&\Omega_{n \,\flat}(f,h) = \\
&\biggl\{ \int_{|\xi|=1} \tr\biggl( \sum \frac{1}{\a!\b!\d!}
 (D^\b_x f)(D^{\a+\d}_x h)\del^{\a+\b}_\xi(\s_L^F)\del^\d_\xi(\s_L^F) \biggr)
 \,d^{n-1}\xi \biggr\} \,d\,x,
\end{align*}
with the sum taken over $|\a|+|\b|+|\d| = n$, $1 \leq |\b|$, $1 \leq |\d|$.
\end{lemma}

To better handle the previous expression, we consider
$$
\phi(\xi,u,v) :=
\sum\frac{u^\b v^{\a+\d}}{\a!\b!\d!}
\tr(\del^{\a+\b}_\xi(\s_L^F)\del^\d_\xi(\s_L^F))
$$
with the sum as before; in this way
$$
\Omega_{n \,\flat}(f,h)
= \Bigl\{\sum A_{a,b}(D_x^a f) (D_x^b h) \Bigr\} \,d \,x
$$
with
$$
\sum A_{a,b} u^a v^b = \int_{|\xi|=1} \phi(\xi,u,v) \,d^{n-1}\xi.
$$
Our next task is to compute the previous integral.
Instead of a direct approach to compute
$\tr\bigl(\del^{\a+\b}_\xi(\s_L^F) \del^\d_\xi(\s_L^F)\bigr)$,
we shall try to use the Taylor expansion of the function
$$
\psi(\xi,\eta) = \tr(\s_L^F(\xi) \s_L^F(\eta)),
$$
on the diagonal point $(\xi,\xi)$, as suggested in \cite{Con1} for the 4-dimensional case. Indeed, with
$\a=(\a_1,\dots,\a_n,\a_{n+1},\dots,\a_{n+n}) = (\b,\d)$, we have
\begin{align*}
\psi(\xi + u, \eta + v)
&= \sum_{|\a|\geq 0} \frac{(u,v)^\a}{\a!}
     \tr \biggl(\del^\a_{(\xi,\eta)}(\s_L^F(\xi) \s_L^F(\eta))
                \biggr|_{\eta=\xi} \biggr)
\\
&= \sum_{|\b|+|\d|\geq 0} \frac{u^\b}{\b!} \frac{v^\d}{\d!}
     \tr \biggl(\del_{(\xi,\eta)}^{(\b,0)}(\s_L^F(\xi))
                \del_{(\xi,\eta)}^{(0,\d)}(\s_L^F(\eta))
                \biggr|_{\eta=\xi} \biggr) \\
&= \sum_{|\b|+|\d|\geq 0} \frac{u^\b}{\b!} \frac{v^\d}{\d!}
     \tr \biggl(\del_\xi^\b(\s_L^F(\xi)) \del_\eta^\d(\s_L^F(\eta))
                \biggr|_{\eta=\xi} \biggr).
\end{align*}

Therefore, the term of order $n$ in the Taylor expansion of $\psi$ at
the point $(\xi,\xi)$ is given by
\begin{align*}
& T_n \psi(\xi,\xi,u,v) =
\sum_{|\b|+|\d|=n} \frac{u^\b}{\b!} \frac{v^\d}{\d!}
     \tr\biggl(\del_\xi^\b(\s_L^F(\xi)) \del_\eta^\d(\s_L^F(\eta))
               \biggr|_{\eta=\xi} \biggr)\\
&= \sum_{{|\b|+|\d|=n}\atop{|\b|\geq 1,|\d|\geq 1}}
               \frac{u^\b}{\b!} \frac{v^\d}{\d!}
     \tr\biggl(\del^\b_\xi(\s_L^F(\xi)) \del_\eta^\d(\s_L^F(\eta))
               \biggr|_{\eta=\xi} \biggr) \\
& \quad + \sum_{|\b|=n} \frac{u^\b}{\b!}
     \tr\bigl(\del^\b(\s_L^F(\xi)) \s_L^F(\xi)\bigr)
  + \sum_{|\d|=n} \frac{v^\d}{\d!}
     \tr\biggl(\s_L^F(\xi) \del_\eta^\d(\s_L^F(\eta))
               \biggr|_{\eta=\xi}\biggr).
\end{align*}

We denote by $T'_n \psi(\xi,\eta,u,v)$ the term of order $n$ in the
Taylor expansion of $\psi(\xi,\eta)$ minus the terms with only powers
of $u$ or only powers of $v$. That is to say,
$$
T'_n \psi(\xi,\eta,u,v) =
\sum_{{|\b|+|\d|=n}\atop{|\b|\geq 1,|\d|\geq 1}}
      \frac{u^\b}{\b!} \frac{v^\d}{\d!}
     \tr\bigl(\del^\b_\xi(\s_L^F(\xi)) \del^\d_\eta(\s_L^F(\eta))\bigr).
$$
Now
\begin{align*}
&T'_n \psi(\xi,\eta,u+v,v)\\
&= \sum_{{{|\b|+|\d| = n}\atop{|\b|\geq 1,|\d|\geq 1}}}
        \frac{(u+v)^\b v^\d}{\b!\d!}
     \tr\bigl(\del^\b_\xi(\s_L^F(\xi)) \del^\d_\xi(\s_L^F(\eta))\bigr) \\
&= \sum_{{|\b|+|\d| = n}\atop{|\b|\geq 1, |\d| \geq 1}}
   \sum_{\b'+\b''=\b} \frac{u^{\b'}v^{\b''+\d}}{\b'!\b''!\d!}
     \tr\bigl(\del^\b_\xi(\s_L^F(\xi)) \del^\d_\eta(\s_L^F(\eta))\bigr) \\
&= \sum_{{|\b|+|\d| = n}\atop{|\b|\geq 1, |\d| \geq 1}}
   \sum_{{\b'+\b''=\b}\atop{\b'\not= 0}} \frac{u^{\b'}v^{\b''+\d}}{\b'!\b''!\d!}
     \tr\bigl(\del^\b_\xi(\s_L^F(\xi)) \del^\d_\eta(\s_L^F(\eta))\bigr) \\
&\qquad + \sum_{{|\b|+|\d|=n}\atop{|\b|\geq 1, |\d|\geq 1}}
          \frac{v^{\b+\d}}{\b!\d!}
     \tr(\del^{\b}_\xi(\s_L^F(\xi))\del^\d_\eta(\s_L^F(\eta))) \\
&= \sum_{{|\b'|+|\b''|+|\d|=n}\atop{|\b''|\geq1,|\d|\geq1}}
          \frac{u^{\b'}v^{\b''+\d}}{\b'!\b''!\d!}
     \tr\bigl(\del^{\b'+\b''}_\xi(\s_L^F(\xi)) \del^\d_\eta(\s_L^F(\eta)) \bigr)
   + T'_n \psi(\xi,\eta,v,v).
\end{align*}

Therefore, by taking $\eta = \xi$ we obtain:
$$
T'_n \psi(\xi,\xi,u+v,v) - T'_n \psi(\xi,\xi,v,v) = \phi(\xi,u,v),
$$
which means that to compute $\phi$ we only need to compute the term of
order $n$ in the Taylor expansion of $\psi$ at $(\xi,\xi)$ and forget
about the terms with only powers of $u$ or $v$. Summarizing, we have

\begin{theorem}
\label{theoremomeganflat}
$$
\Omega_{n \,\flat}(f,h)
= \Bigl(\sum A_{a,b} (D_x^a f)(D_x^b h) \Bigr) \,d\,x,
$$
where
$$
\sum A_{a,b} u^a v^b =
\int_{|\xi|=1} \bigl(T'_n \psi(\xi,\xi,u+v,v) -
                     T'_n \psi(\xi,\xi,v,v) \bigr) \,d^{n-1}\xi
$$
and $T'_n \psi(\xi,\eta,u,v)$ is the term of order $n$ in the Taylor
expansion of the function
$\psi(\xi,\eta) = \tr(\s^F_L(\xi) \s^F_L(\eta))$ without the terms
with only powers of $u$ or only powers of $v$, at the
point $(\xi,\eta)$.
\end{theorem}

%%%%%%%%%%%%%%%%%%%%%%%%%%%%%%%%%%%%%%%%%%%%%%%%%%%%%%%%%%%%%%%%%%%

%% subsection: A general formula for $\tr(\s_L^F(\xi)\s_L^F(\eta))$ %%

\noindent{\bf A general formula for $\tr(\s_L^F(\xi)\s_L^F(\eta))$}

Because of Theorem~\nnn{theoremomeganflat}, to obtain an explicit expression for $\Omega_n$, at least for the flat
case, it is necessary to study $\tr(\s_L^F(\xi) \s_L^F(\eta))$ for $\xi$ and $\eta$ not zero in $T^*_xM$.
Note that, ignoring the space of harmonic forms at each level, the operator $F=(d\d-\d d)\Delta^{-1}$, and
its symbol, are defined and preserve the subspace of $m$-forms for any $0 \leq m \leq n$.

\begin{theorem}
\label{theoremtraces0}
With $\s^F_L(\xi) \s^F_L(\eta)$ acting on $m$-forms we have:
\begin{equation}
\label{trsxiseta}
\tr(\s^F_L(\xi) \s^F_L(\eta)) =
a_{n,m} \frac{\<\xi,\eta>^2}{|\xi|^2|\eta|^2} + b_{n,m},
\end{equation}
where $\<\xi,\eta>$ represents the inner product $g^m(\xi,\eta)$
given by the metric and
$$
\binom{n}{m} - a_{n,m} = b_{n,m} =
\binom{n-2}{m-2} +\binom{n-2}{m} - 2 \binom{n-2}{m-1}.
$$
\end{theorem}

To prove it, we need to set some notation and preliminary results. In
the following $\xi_i$, $\eta_i$, $\xi$ and $\eta$ belong to $T^*_x M$
and are not zero.  Also we will represent the trace operator acting on the different
subspaces with the same symbol, making clear out of the context on which dimension is it acting.

By using the trace identity $\tr(A B)=\tr(B A)$ and the relation
\begin{equation}
\label{innerxieta}
\eps_{m-1}(\xi)\iota_m(\eta) + \iota_{m+1}(\eta)\eps_m(\xi)
= \<\xi,\eta> I_m
\end{equation}
with $I_m$ the identity on $m$-forms, we know, for every $m\geq 1$,
$$
\tr(\iota_m(\eta) \eps_{m-1}(\xi))
+ \tr(\iota_{m+1}(\eta) \eps_m(\xi))
=  \binom{n}{m} \, \<\xi,\eta>.
$$
In the particular case $m = 0$, we have
$\iota_1(\eta) \eps_0(\xi) = \<\xi,\eta>I_0$, so
$\tr(\iota_1(\eta) \eps_0(\xi)) = \<\xi,\eta>$, and hence

\begin{lemma}
\label{trie}
For every $m \geq 0$,
\begin{align*}
\tr(\iota_{m+1}(\eta) \eps_{m}(\xi))
&= \<\xi,\eta> \biggl(\binom{n}{m} - \binom{n}{m-1} +\cdots+
   (-1)^m \binom{n}{0} \biggr) \\
&=  A_{n,m} \, \<\xi,\eta>.
\end{align*}
\end{lemma}

By Proposition~\ref{leadingsymbolofF},
$\s_L^F(x,\xi) = |\xi|^{-2}(\eps_{m-1}(\xi)\iota_{m}(\xi)
- \iota_{m+1}(\xi)\eps_{m}(\xi))$, so
\begin{align*}
\tr(\s_L^F(x,\xi))
&= |\xi|^{-2} \bigl(\tr(\iota_m(\xi)\eps_{m-1}(\xi))
-\tr(\iota_{m+1}(\xi)\eps_{m}(\xi))\bigr) \\
&= |\xi|^{-2} \biggl(\binom{n}{m} |\xi|^2
- 2\,\tr(\iota_{m+1}(\xi)\eps_{m}(\xi))\biggr).
\end{align*}
By \nnn{innerxieta},
\begin{equation}
\label{ei-ie=2eiI}
\eps_{m-1}(\xi)\iota_m(\xi) -\iota_{m+1}(\xi)\eps_m(\xi) =
2 \eps_{m-1}(\xi)\iota_m(\xi) - |\xi|^2 I_m,
\end{equation}
hence Proposition~\ref{leadingsymbolofF} implies
\begin{align*}
\tr(\s_L^F(x,\xi))
&= |\xi|^{-2}\tr(\eps_{m-1}(\xi)\iota_m(\xi)
    - \iota_{m+1}(\xi)\eps_m(\xi)) \\
&= 2|\xi|^{-2}\tr(\eps_{m-1}(\xi)\iota_m(\xi)) - \binom{n}{m}
= 2 A_{n,m-1} - \binom{n}{m},
\end{align*}
and we obtain as a consequence of the previous Lemma:

\begin{corollary}
For every $(x,\xi) \in T^*M$,
$$
\tr(\s_L(F)(x,\xi)) = 2A_{n,m-1} - \binom{n}{m}.
$$
\end{corollary}

\begin{proof}[Proof of Theorem~\ref{theoremtraces0}]
Using \nnn{ei-ie=2eiI}, and Lemma~\ref{trie} we deduce
\begin{align*}
&\tr\bigl(\s_L(F)(x,\xi) \s_L(F)(x,\eta)\bigr) \\
&= |\xi|^{-2} |\eta|^{-2}
   \tr\bigl((2 \eps_{m-1}(\xi)\iota_m(\xi) - |\xi|^2 I_m)
            (2 \eps_{m-1}(\eta)\iota_m(\eta) - |\eta|^2 I_m)\bigr)\\
&= 4|\xi|^{-2} |\eta|^{-2}
   \tr(\eps_{m-1}(\xi)\iota_m(\xi)\eps_{m-1}(\eta)\iota_m(\eta))
   - 4A_{n,m-1} + \binom{n}{m}.
\end{align*}

Considering for any $m \geq 1$ the quantity
$$
a_m(\xi_1,\xi_2,\eta_1,\eta_2) =
\tr(\eps_{m-1}(\xi_1)\iota_m(\xi_2)\eps_{m-1}(\eta_1)\iota_m(\eta_2)),
$$
we have the relation
\begin{align*}
& a_m(\xi_1,\xi_2,\eta_1,\eta_2) \\
& = \tr\bigl((\<\xi_1,\xi_2>I_m - \iota_{m+1}(\xi_2)\eps_m(\xi_1))
      (\<\eta_1,\eta_2>I_m - \iota_{m+1}(\eta_2)\eps_m(\eta_1))\bigr)
\\
& = \<\xi_1,\xi_2> \<\eta_1,\eta_2> \biggl(\binom{n}{m} -2 A_{n,m}\biggr)
    + \tr(\iota_{m+1}(\xi_2)\eps_{m}(\xi_1)
          \iota_{m+1}(\eta_2)\eps_{m}(\eta_1)),
\end{align*}
which implies
$$
a_{m+1}(\eta_1,\xi_2,\xi_1,\eta_2) =
a_{m}(\xi_1,\xi_2,\eta_1,\eta_2) +
\<\xi_1,\xi_2> \<\eta_1,\eta_2> \biggl(2 A_{n,m} - \binom{n}{m} \biggr),
$$
where
\begin{align*}
a_1(\xi_1,\xi_2,\eta_1,\eta_2)
&= \tr(\eps_0(\xi_1)\iota_1(\xi_2) \eps_0(\eta_1)\iota_1(\eta_2)) \\
&= \tr(\iota_1(\eta_2)\eps_0(\xi_1) \iota_1(\xi_2)\eps_0(\eta_1)) \\
&= \tr(\<\eta_2,\xi_1> \<\xi_2,\eta_1> I_0)
= \<\eta_2,\xi_1> \<\xi_2,\eta_1>.
\end{align*}

Thus, because of the relation $$
\tr(\s_L(F)(x,\xi)\s_L(F)(x,\eta)) = 4
a_m(\xi,\xi,\eta,\eta)|\xi|^{-2}|\eta|^{-2} + \binom{n}{m} - 4
A_{n,m-1}, $$ we have a recursive way of computing the left hand
side. This is enough to prove that
$\tr(\s^F_L(\xi)\s^F_L(\eta)=a_{n,m}\<\xi,\eta^2>|\xi|^{-2}|\eta|^{-2}+b_{n,m}$
for some constants $a_{n,m}$ and $b_{n,m}$.

Now in the particular case in which $\xi = \eta$ is a member of an
orthonormal basis of 1-forms, the operator under the trace equals
$I_m$, the identity on $m$-forms because of the equality
$\s_L(F)(\xi)\s_L(F)(\xi)=(\eps(\xi)\iota(\xi)
-\iota(\xi)\eps(\xi))^2$ and the following computations. For a
basic $m$-form $e_{j_1} \wedge\cdots\wedge e_{j_m}$, one has $$
\bigl(\eps_{m-1}(\xi)\iota_m(\xi) -
\iota_{m+1}(\xi)\eps_m(\xi)\bigr) (e_{j_1} \wedge\cdots\wedge
e_{j_m}) = e_{j_1} \wedge\cdots\wedge e_{j_m} $$ if $\xi =
e_{j_i}$ for some $i;$  and $$ \bigl(\eps_{m-1}(\xi)\iota_m(\xi) -
\iota_{m+1}(\xi)\eps_m(\xi)\bigr) (e_{j_1} \wedge\cdots\wedge
e_{j_m}) = - e_{j_1} \wedge\cdots\wedge e_{j_m} $$ if $\xi \neq
e_{j_i}$ for every $i$.  Therefore $a_{n,m}+b_{n,m}=\binom{n}{m}$.

If both $\xi$ and $\eta$ are different members of an orthonormal basis, the term $\<\xi,\eta>$
vanishes and the expression \nnn{trsxiseta} reduces to $b_{n,m}$, in this case,
$$
\s_L(F)(\xi) \s_L(F)(\eta) e_{j_1} \wedge\cdots\wedge e_{j_m} =
e_{j_1} \wedge\cdots\wedge e_{j_m}
$$
if both $\xi,\eta \in \{e_{j_1}, \dots, e_{j_m}\}$ or if both
$\xi,\eta \notin \{e_{j_1}, \dots, e_{j_m}\}$; and
$$
\s_L(F)(\xi) \s_L(F)(\eta) e_{j_1} \wedge\cdots\wedge e_{j_m} =
- e_{j_1} \wedge\cdots\wedge e_{j_m}
$$
if only one of $\xi$, $\eta$ belongs to $\{e_{j_1},\dots,e_{j_m}\}$
and the other does not.

The number of basic $m$-forms containing both $\xi$ and $\eta$ as factors is $\binom{n-2}{m-2}$.
The number of basic $m$-forms containing neither $\xi$ nor $\eta$ is $\binom{n-2}{m}$, and
the number of basic $m$-forms containing exactly one of $\xi$ or $\eta$ as a factor is
$2 \binom{n-2}{m-1}$.
In this way, the value of $b_{n,m}$ is given by the trace of the above operator, which equals
$b_{n,m} = \binom{n-2}{m-2} + \binom{n-2}{m} - 2 \binom{n-2}{m-1}$,
and hence the value for $a_{n,m}$ is
$a_{n,m} = \binom{n}{m} - \binom{n-2}{m-2} - \binom{n-2}{m} +
2 \binom{n-2}{m-1}$,
and the proof is complete.
\end{proof}

We can now restate Theorem~\ref{theoremomeganflat} as

\begin{theorem}
\label{omegaflatcase}
$$
\Omega_{n \,\flat}(f,h)
= \Bigl(\sum A_{a,b} (D_x^a f)(D_x^b h) \Bigr) \,d\,x,
$$
where
$$
\sum A_{a,b} u^a v^b =
\int_{|\xi|=1} \bigl(T'_n\psi(\xi,\xi,u+v,v) -
                     T'_n\psi(\xi,\xi,v,v) \bigr) \,d^{n-1}\xi
$$
and $T'_n \psi(\xi,\eta,u,v)$ is the term of order $n$ in the Taylor
expansion of the function
$\psi(\xi,\eta) = a_{n,l} \<\xi,\eta>^2|\xi|^{-2}|\eta|^{-2} + b_{n,l}$ without the terms
with only powers of $u$ or only powers of $v$, at the point $(\xi,\eta)$.
\end{theorem}

%%%%%%%%%% section: An approach to the general case %%%%%%%%%%%

\section{An approach to the general case}
\label{generalapproach}

In this section we want to address a possible approach to the general case.  This approach is the same used in
\cite{Con1} to compute $\Omega_4$ and in \cite{Irazu} for $\Omega_6$.  The first step to obtain
$\Omega_n$ is to compute it explicitly in a flat metric. That is precisely the result given by
Theorem~\ref{omegaflatcase}.  Then, by changing the metric conformally, one
obtains the expression for $\Omega_n$ in the conformally flat case as follows.
If $g$ is a conformally flat metric, then $g=e^{2 \eta}g_{\flat}$ for a $C^\infty$ function $\eta$ and
$g_\flat$ the flat metric on $M.$  One expresses $\Omega_{n,\flat}(f,h)$ in terms of the new metric $g$,
and obtains an expression in terms of the covariant derivatives of $\eta,$ $f,$ $h,$ and the Riemann curvature
tensor $R.$  By using the conformal change equation for the Ricci tensor:
\begin{equation}
\label{invariantize} \llcovd{\eta}{i}{j} = - V_{ij} -
\lcovd{\eta}{i}\,\lcovd{\eta}{j} +
\frac{1}{2}\lcovd{\eta}{k}\,\ucovd{\eta}{k}\,g_{ij},
\end{equation}
one replaces
the second covariant derivatives on $\eta$ with terms with the Ricci tensor.
The result is the expression for $\Omega_n(f,h)$ in the conformally flat metric $g.$

In \nnn{invariantize}
$V$ represents a normalized translation of the Ricci tensor,
useful in conformal geometry, given in terms of the normalized scalar curvature $J$ by
$$
V = \frac{\Rc - J g}{n - 2}
\quad \hbox{with}\quad
J= \frac{\Sc}{2(n - 1)},
$$
and the indices after the semicolon represents covariant derivatives,
$\llcovd{\eta}{i}{j}=\Na_j \Na_i \eta.$

To study the conformal invariance,
there is no need to study the whole conformal deformation. It is
enough to study the conformal deformation up to order one in $\eta$
as follows. If we set $\hat g = e^{2 z \eta}g_\flat$ where
$\eta \in C^\infty(M)$ and $z$ a constant, then the conformal
variation of each expression is a polynomial in~$z$
whose coefficients are expressions in the metric and the conformal factor
$\eta$ (actually, this is an abuse of the language since the
conformal factor is $e^{2 z \eta}$). In this way, the conformal
deformation up to order one in $\eta$ is given by
$\frac{d}{d z}\bigr|_{z=0}.$   If the conformal deformation of a natural tensor or a differential
operator is zero up to order one then, by integration it is fully invariant, for details see \cite{Bra4}.

The last step is to take the expression
obtained in the conformally flat case and consider it in a general metric.  Finally, study its conformal
variation (up to order one), and then look for those terms that must be added
to obtain a conformally invariant expression for a general conformally curved metric.

Following the spirit of \cite{Con1}, one reorder
the summation of the expression given by \nnn{wresf0f1f2} for $\Omega_n(f,h)$ by looking at
$j+k+|\a'|$ in $\set{0,1,\dots,n-2}$. In this way $\tr(\s_{-n}([F,f][F,h]))$ can be written as
$\sum_{p=0}^{n-2} \tr(\s^{(p)}_{-n})$,
where
\begin{align}
&\tr(\s^{(p)}_{-n}) = \nonumber \\
&\sum \frac{1}{\a'!\a''!\b!\d!} D^\b_x(f) D_x^{\a''+\d}(h)
     \tr\bigl(\del_\xi^{\a'+\a''+\b}(\s^F_{-j})
              \del_\xi^\d (D_x^{\a'}(\s^F_{-k})) \bigr),
\label{s_np}
\end{align}
with the sum taken over $j+k+|\a'| = p$.

Each $\s_{-j}^F$ is a $\binom{n}{l} \times \binom{n}{l}$
matrix with $\s^F_0$ only invoking $g_{i j}(x)$ and the $\s^F_{-j}$ for
$j\geq1$ are polynomial expressions on the partial derivatives of the metric at
$x$, see Lemma~\nnn{lemma1'}.

By \nnn{s_np}, $\tr(\s^{(0)}_{-n})$ is given in terms of
$g_{i j}(x)$ (without any of its partial derivatives) since only contains terms of the
form $\tr(\del^{\a''+\b}_\xi (\s_L^F)\del^{\d}_{\xi}(\s_L^F))$.  Also by
\nnn{s_np}, any $\tr(\s^{(p)}_{-n})$ for $p \geq 1$ is a
polynomial expression on the partial derivatives of the metric at $x$ with coefficients
depending smoothly on $g_{ij}(x)$, indeed, each of these $\tr(\s^{(p)}_{-n})$ only
invoke terms of the form
$\tr(\del^{\a'+\a''+\b}_\xi(\s_{-j}^F) \del^{\d}_{\xi}(D^{\a'}_x(\s_{-k}^F)))$

Evidently, these properties are preserved after integration over the
variable $\xi$. By choosing the coordinates $x^j$ to be geodesic
normal coordinates at the point $x$, one can assume that
$g_{i j}(x) = \d^j_i$, that the first partial derivatives of the metric vanish at
$x$, and that the higher partial derivatives are expressed in term of the Riemannian
curvature and its covariant derivatives \cite{Gray}. In this way, we have generalized
Lemma 4 in \cite{Con1} to any even dimension by splitting $\Omega_n(f,h)$ as the sum of
$$
\int_{|\xi|=1} \tr(\s^{(0)}_{-n})\,d^{n-1}\xi\,d\,x= B_n(\Na^a d f, \Na^b d h)
$$
and
$$
\left\{\int_{|\xi|=1} \tr(\s^{(1)}_{-n} + \cdots +
\s^{(n-2)}_{-n})\,d^{n-1}\xi\right\}\,d \,x
= C_n(R,d f,d h).
$$
In this way one has
\begin{lemma}
\label{lemmaBC}
There exists a universal bilinear expression $B_n(\Na^\a d f, \Na^\b d h)$
for some multi-indices $\a$ and $\b$ and an expression $C_n(R,d f,d h)$ linear in $f$ and
$h$ such that
\begin{equation}
\label{OmegaBC}
\Omega_n(f,h) = \bigl\{B_n(\Na^\a d f,\Na^\b d h) + C_n(R,d f,d h)\bigr\} \,d \,x
\end{equation}
where $R$ is the Riemannian curvature tensor, $\Na$ the covariant
differentiation, and $|\a| + |\b| = n - 2$.
\end{lemma}

Note that covariant derivatives do not commute, making the notation
$\Na^\a d f$ somehow ambiguous. It is assumed that $C(R,d f,d h)$ will
absorb any ambiguity.  Also, because of the restriction $|\b|\geq1$, $|\d|\geq1$ in
Lemma~\nnn{lemmas_nSS}, $\Omega_n(f,h)$ depends only on $df$ and $dh$ rather than on $f$
and $h$, that is why the right hand side of \nnn{OmegaBC} is given in terms of $d f$ and
$dh$.

In the particular case of a flat manifold, all the terms involving
curvature vanish, in that situation $C_n(R,d f,d h)$ reduces to zero and
the covariant derivatives do commute. In particular,
$B_n(\Na^\a d f, \Na^\b d h)$ equals $\Omega_n(f,h)$.

\begin{lemma}
If $M$ is a flat manifold,
$$
\Omega_n(f,h) = B_n(\Na^\a d f, \Na^\b d h) \,d \,x.
$$
\end{lemma}

{\bf How to compute $C_n(R,d f,d h)$?}
The expression $\omega_n(f,h)$ in \nnn{wresf0f1f2} is a sum of homogeneous polynomials in the
ingredients $\Na^\a d f$, $\Na^\b d h$, and $\Na^\ga R$ for multi-indices $\a,\b$, and $\ga$,
in the following sense, each monomial must satisfies the {\bf homogeneity condition} given by the
rule (see \cite{Bra4}):
$$
2 k_R + k_\Na = n
$$
where $k_R$ denotes its degree in~$R$ and $k_\Na$ its degree in~$\Na$.
Also, for covariant derivatives we count all of the derivatives on
$R$, $f$, and $h$, and any occurrence of $Rc$, or $Sc$ as an occurrence of~$R$.
By closing under addition, we denote by $\cP_n$ the space of these
polynomials.

Because $|\b|\geq 1$, and $|\d|\geq 1$ in Lemma~\ref{lemmas_nSS}, we
have $k_\Na \geq 2$ and hence $k_R \leq (n-2)/2$.
We say that $Q$ is in $\cP_{n,l}$ if $Q$ can be written as a sum of
monomials with $k_R\geq l$, or equivalently, $k_\Na \leq n - 2l$. We
have
$$
\cP_n = \cP_{n,0} \supseteq \cP_{n,1} \supseteq \cP_{n,2}
\supseteq\cdots\supseteq \cP_{n,\tfrac{n-2}{2}},
$$
and $\cP_{n,l} = 0$ for $l > (n-2)/2$.

There is an important observation to make. An expression which a
priori appears to be in, say $\cP_{6,1}$, may actually be in a
subspace of it, like $\cP_{6,2}$; for example,
$$
\underbrace{\lcovd{f}{i} \lllcovd{h}{j}{k}{l} W^{ijlk}}_{\in\,\cP_{6,1}} =
\underbrace{\lcovd{f}{i} \lcovd{h}{j} V_{kl} W^{ikjl}
+\lcovd{f}{i}\lcovd{h}{j} W^i{}_{klm} W^j{}_{klm}}_
{\in \, \cP_{6,2}},
$$
by reordering covariant derivatives and making use of the symmetries of the
Weyl tensor $W$.
Because of this filtration, it is necessary to adopt a fixed
convention on how the indices should be placed when representing
each expression in its index notation.
For example, $\lllcovd{f}{i}{j}{k} \uuucovd{h}{i}{j}{k}$ will be
preferred over $\lllcovd{f}{i}{k}{j} \uuucovd{h}{i}{j}{k}$ or any
other variation.  Lemma~\nnn{lemmaBC} is restated as

\begin{lemma}
There exists a universal bilinear form $Q_n(df,dh)$ in
$\cP_{n,0} \smallsetminus \cP_{n,1}$ and a form $Q_{R,n}(df,dh)$ in
$\cP_{n,1}$ such that
$$
\omega_n(f,h) = Q_n(df,dh) + Q_{R,n}(df,dh)
$$
where $R$ is the Riemannian curvature tensor.
\end{lemma}

Once one has computed the expression for $\Omega_n$ in a flat metric $g$, one
express it in a conformally related metric $\hat g$.  By doing that, all of the terms
for $Q_{R,n}$ which do not involve a factor of the Weyl tensor $W$ are computed.
In this way, up to the conformally flat case, one obtains the expression for $\Omega_n$.
Next, one uses the homogeneity condition $2 k_\Na+k_R=n$ to list the possible
terms to be considered in the general conformally curved case.  These terms must be such that they contain
the Weyl tensor as a factor, and
their conformal variation compensate that of the already known terms for $\Omega_n$.

In the particular case $n=4$, as found in \cite{Con1}, the homogeneity condition shows that
no extra terms exist and the expression found in the conformally flat case suffices.
Indeed, $\Omega_4(f,h)$ is given by a constant multiple of
$$
-\llucovd{f}{i}{j}{j}\ucovd{h}{i} - \frac{1}{2}\lucovd{f}{i}{i}\lucovd{h}{j}{j}
- \llcovd{f}{i}{j}\uucovd{h}{i}{j} - \lcovd{f}{i}\ulucovd{h}{i}{j}{j}
+2 \lcovd{f}{i}\ucovd{h}{i}J,
$$
where $2(n-1)J=Sc$.  The only terms to be added, $\lcovd{f}{i} \ucovd{h}{i} W_*{}^*{}_*{}^*$ or
$\lcovd{f}{i}\lcovd{h}{j}W^i{}_*{}^j{}^*$, are zero because of the symmetries of the Weyl tensor $W$ (here each $*$ represents an index).

In the case $n=6$, as in \cite[Section~6]{Irazu}, six extra terms
are needed as a result of the homogeneity condition. As shown in
that reference, symmetry, conformal invariance, and even the
Hochschild 2-cocycle property of $\int_M f_0 \Omega_n(f_1,f_2)$
are not enough to compute the coefficients in front of those extra
terms containing the Weyl tensor as a factor, indeed, it is not
possible to determine, out of those three properties,  the
coefficients for $\lcovd{f}{i}\ucovd{h}{j}W_{jklm}W^{jklm}$ and
$\lcovd{f}{i}\lcovd{h}{j}W^i{}{}_{klm}W^{jklm}$.

%%%%%%%%%%%%%%%%%%%%

\appendix
\section{$\s(F)$ and the partial derivatives of the metric}

One of the basis of the result of \cite{Con1} is its Lemma~3,
where $M$ is a 4-dimensional manifold and the operator $F$, as
well as its symbol, are acting on $2$-forms. Our objective in this
paragraph is to give a more general version of this lemma, in
which we try to understand the behavior up to any sub-index, of
the component $\s_{-k}$ of the total symbol of $F.$  At the same
time, we give a detailed computation of the total symbol of
$\Delta^{-1}$ in terms of the total symbol of $\Delta$, and an
explicit recursive computation of the total symbol of $F$ in terms
of the total symbols of $\Delta^{-1}$ and of $d \d - \d d$ in some
given local coordinates.

To simplify the typing, we will abbreviate using
the notation $T \sim \cL_k$ to mean ``\emph{the term $T$ has the property of
being linear in the $k$-partial derivatives of the metric at $x$ with coefficients
depending smoothly on the $g_{ij}(x)$}''.  By $T\sim (\cL_k)^r$ we mean
`` \emph{the term $T$ has the property of being a sum of products of $r$ linear
expressions of the type $\cL_k$}''.  By $T \sim \cL_0$ we mean
``\emph{the term $T$ only invokes $g_{ij}(x)$}'', thus if
$T\sim (\cL_0)(\cL_k)$ then $T \sim \cL_k$.

The first step is an recursive computation of the symbol of $\Delta^{-1}$
in terms of the symbol of $\Delta$.  We use the following reasoning, justified by
\cite[Corollary~1.4.3]{Treves}.  From \cite[Lemma~2.4.4]{Gilkey}, the
operator $\Delta = d \d + \d d$ acting on $\La^{l} T^*M$ can
be expressed in a given system of local coordinates in the form
$\s(d\d + \d d) = p_2 + p_1 + p_0$,
where each $p_i$ is a $\binom{n}{l} \times \binom{n}{l}$ matrix such that $p_2 \sim \cL_0$,
$p_1\sim  \cL_1$, and $p_0 \sim \cL_2 + (\cL_1)^2$.

Let
$$
\s(\Delta^{-1}) = \sum_{k\leq-2} \s_k(\Delta^{-1}) = \sum_{k \leq-2} r_k
$$
be the total symbol of $\Delta^{-1}$, hence
\begin{lemma}
\label{lemasymdelta-1}
The total symbol of $\Delta^{-1}$ is given by the recursive formula
$r_{-2}=p_2^{-1}$,
$r_{-3}= p_2^{-1}(-p_1 p_2^{-1} + \del_\xi^1p_2 \,D_x^1(p_2^{-1}))$,
and in general, for all $k \geq 2$,
\begin{align*}
r_{-(k+2)} &= p_2^{-1} \bigl(-p_1 r_{-(k+1)} - p_0 r_{-k} \\
&\qquad   - (\del^1_\xi p_2 \,D^1_x r_{-(k+1)} + \del_\xi^1 p_1 \,D_x^1 r_{-k})
       - \tfrac12 \del_\xi^2 p_2 \,D_x^2 r_{-k} \bigr).
\end{align*}
\end{lemma}

\begin{proof}
\begin{align}
1&= \sum_\a \frac{1}{\a!} \del_\xi^\a \s(\Delta) \, D_x^\a \s(\Delta^{-1})
  = \sum_{|\a|=0}^2 \frac{1}{\a!}
        \del_\xi^\a \s(\Delta) \, D_x^\a \s(\Delta^{-1})
\nonumber \\
&= (p_2 + p_1 + p_0)(r_{-2} + r_{-3} + r_{-4} + \cdots)
\nonumber \\
&\qquad + (\del_\xi^1 p_2 + \del_\xi^1 p_1)
         (D_x^1 r_{-2} + D_x^1 r_{-3}+ D_x^1 r_{-4} + \cdots)
\nonumber \\
&\qquad + \tfrac12 (\del_\xi^2 p_2)
         (D_x^2 r_{-2} + D_x^2 r_{-3} + D_x^2 r_{-4} + \cdots)
\label{sDeltaInvDelta} \\
&= [p_2r_{-2}] + [p_2r_{-3} + p_1r_{-2} + \del_\xi^1 p_2\,D_x^1 r_{-2}]
\nonumber \\
&\qquad + \sum_{k=2} \bigl[ p_2r_{-(k+2)} + p_1r_{-(k+1)} + p_0r_{-k}
\nonumber \\
&\qquad\quad + (\del_\xi^1 p_2\,D_x^1 r_{-(k+1)}
              + \del_\xi^1 p_1\,D_x^1 r_{-k})
             + \tfrac12 \del_\xi^2 p_2\,D_x^2 r_{-k} \bigr] + \cdots
\nonumber
\end{align}
where $D_\b^m = \sum_{|\a|=m} D_\b^\a$ and each term in the summation
has homogeneous order $k$ in the variable~$\xi$.
The lemma follows by solving for $r_{-k}$ at each level of the decomposition.
\end{proof}

\begin{lemma}
\label{lemma1'} Acting on middle dimension forms, the total symbol
of $\Delta^{-1}$ is given as a sum of $\binom{n}{l} \times
\binom{n}{l}$ matrices of the form $\s(\Delta^{-1}) = r_{-2} +
r_{-3} + r_{-4} + \cdots$ where $r_{-2} \sim \cL_0$, $r_{-3} \sim
\cL_1$, $r_{-4} \sim \cL_2 + (\cL_1)^2$, and in general, for
$k\geq 0$, $$ r_{-(k+2)} \sim \sum (\cL_k)^{a_k}
(\cL_{k-1})^{a_{k-1}} \cdots (\cL_1)^{a_1} $$ with the summation
taken over $0 \leq a_i \leq k$, for every $i$ and $a_1 + 2a_2 +
\cdots + k a_k = k$.
\end{lemma}

\begin{proof}
The differentiation $\del_\xi^\a$ does not alter the properties of
$p_k$, so that, for instance, $\del_\xi^\a p_1 \sim \cL_1$. The
differentiation $D_x^\a$ behaves in the following way
$D_x^\a(\cL_k) \sim \cL_{k+|\a|}$ and
$D_x(\cL_k \cL_j) \sim \cL_{k+1}\cL_j + \cL_k\cL_{j+1}$.   The lemma follows
from the next inductive reasoning:

The first three cases read directly from the
expressions for $r_{-j}$ in Lemma~\ref{lemasymdelta-1}, with $j = 2,3,4$. For the general case, assume
\begin{align}
r_{-(k+2)}
&\sim \cL_1 r_{-(k+1)} + \cL_2 r_{-k}\,D^1_x r_{-(k+1)} +
\cL_1 \,D_x^1 r_{-k} + D_x^2 r_{-k}
\nonumber \\
&\sim \cL_1 \sum (\cL_{k-1})^{a_{k-1}}\cdots (\cL_1)^{a_1}
    + \cL_2 \sum (\cL_{k-2})^{b_{k-2}}\cdots (\cL_1)^{b_1}
\nonumber \\
&\qquad + D_x^1 \sum (\cL_{k-1})^{a_{k-1}}\cdots (\cL_1)^{a_1}
        + \cL_1\,D_x^1 \sum (\cL_{k-2})^{b_{k-2}}\cdots (\cL_1)^{b_1}
\nonumber \\
&\qquad + D_x^2 \sum (\cL_{k-2})^{b_{k-2}}\cdots (\cL_1)^{b_1}
\label{r-(k+2)}
\end{align}
where $a_1 + 2  a_2 + \cdots + (k-1)a_{k-1} = k - 1$ and
$b_1 + 2 b_2 + \cdots + (k-2)b_{k-2} = k - 2$.

First note that
$$
\cL_1 \sum (\cL_{k-1})^{a_{k-1}} \cdots (\cL_1)^{a_1}
= \sum (\cL_{k-1})^{a_{k-1}} \cdots (\cL_1)^{1+a_1}
$$
where $(1+a_1) + 2\a_2 + \cdots + (k-1)a_{k-1} = 1 + k-1 = k$, and
$$
\cL_2 \sum (\cL_{k-2})^{b_{k-2}} \cdots (\cL_1)^{b_1}
= \sum (\cL_{k-2})^{b_{k-2}} \cdots (\cL_2)^{1+b_2} (\cL_1)^{b_1}
$$
where $b_1 + 2(1+b_2) + \cdots + (k-2)b_{k-2} = 2 + k-2 = k$.
In this way, the first two terms of \nnn{r-(k+2)} have the required property.

Next,
\begin{align}
& D_x^1 \sum (\cL_{k-1})^{a_{k-1}} \cdots (\cL_1)^{a_1}
\nonumber \\
&= \sum \cL_k (\cL_{k-1})^{a_{k-1}} \cdots (\cL_1)^{a_1} +
   \sum (\cL_{k-1})^{a_{k-1}+1} (\cL_{k-2})^{a_{k-2}-1} \cdots
        (\cL_1)^{a_1} + \cdots
\nonumber \\
&\qquad + \sum (\cL_{k-1})^{a_{k-1}} \cdots
               (\cL_2)^{a_2+1} (\cL_1)^{a_1-1},
\label{Doftypes}
\end{align}
where each of the summations is taken over a set of the form
$(d_1,\dots,d_k)$ with $d_1 + 2 d_2 + \cdots + kd_k = k$.
Indeed,
\begin{gather*}
a_1 + \cdots + (k-1)(a_{k-1}-1) + k = k-1 - (k-1) + k = k, \\
a_1 + \cdots + (k-2)(a_{k-2}-1) + (k-1)(a_{k-1}+1) =
k-1 - (k-2) + (k-1) = k, \\
\intertext{and}
(a_1 - 1) + 2(a_2 + 1) + \cdots + (k-1)a_{k-1} = 1 + k-1 = k.
\end{gather*}

The previous reasoning shows that
$D_x^1 \Bigl(\sum (\cL_{k-2})^{b_{k-2}} \cdots (\cL_1)^{b_1}\Bigr)$
is a sum over a set of indices of the form
$(a_1,\dots,a_{k-1})$ with $a_1 + 2 a_2 +\cdots+ (k-1)a_{k-1} = k-1$,
thus
$$
\cL_1 \,D_x^1 \sum (\cL_{k-2})^{b_{k-2}} \cdots (\cL_1)^{b_1}
\sim \cL_1 \sum (\cL_{k-1})^{a_{k-1}} \cdots (\cL_1)^{a_1},
$$
which satisfies the desired property as before. For the last term in \nnn{r-(k+2)},
$$
D_x^2 \biggl(\sum (\cL_{k-2})^{b_{k-2}}\cdots (\cL_1)^{b_1}\biggr)
\sim
D_x \biggl( D_x \biggl( \sum (\cL_{k-2})^{b_{k-2}} \cdots (\cL_1)^{b_1}
\biggr) \biggr),
$$
which is of the desired type by applying twice the reasoning
on~\nnn{Doftypes}.
\end{proof}

\begin{lemma}
Acting on middle dimension forms, the total symbol $\s^F$ of $F$,
up to order $-k$ inclusive,  is a $\binom{n}{l} \times \binom{n}{l}$ matrix of the form
$\s^F = \s_0^F + \s_{-1}^F + \cdots + \s_{-k}^F + \cdots$
where $\s_0^F \sim \cL_0$, $\s_{-1}^F\sim \cL_1$,
$\s^F_{-2}\sim \cL_2 + (\cL_1)^2$, and in general
$$
\s_{-k}^F \sim
\sum(\cL_k)^{\a_k}(\cL_{k-1})^{\a_{k-1}} \cdots (\cL_1)^{\a_1}
$$
with the summation taken over $0 \leq \a_i \leq k$, for every $i$ and
$\a_1 + 2 \a_2 + \cdots + k \a_k = k$.
\end{lemma}

\begin{proof}
The operator $d \d - \d d$ satisfies $\s(d \d - \d d) = q_2 + q_1 + q_0$,
with $q_2$ of type $\cL_0$, $q_1$ of type $\cL_1$, and
$q_0$ of type $\cL_2 + (\cL_1)^2$; these results follow exactly as the
result for $\Delta = d\d + \d d$.

As in \nnn{sDeltaInvDelta} we have, for the symbol of
$F = (d \d - \d d)\Delta^{-1}$, the expansion
\begin{align*}
\s(F)
&= [q_2 r_{-2}] + [q_2 r_{-3} + q_1r_{-2} + \del_\xi^1 q_2 \,D_x^1 r_{-2}] \\
&\quad + \sum_{k=2} [q_2 r_{-(k+2)} + q_1 r_{-(k+1)} + q_0 r_{-k} \\
&\qquad + (\del_\xi^1 q_2 \,D_x^1 r_{-(k+1)} + \del_\xi^1 q_1 \,D_x^1 r_{-k})
        + \tfrac12 \del_\xi^2 q_2 \,D_x^2 r_{-k}] + \cdots \\
&= \s^F_0 + \s^F_{-1} +\s^F_{-2} + \cdots .
\end{align*}
{}From this last equality and the properties of $q_i$ and $r_j$, we get that
$\s_0^F \sim \cL_0$, $\s_{-1}^F \sim \cL_1$,
$\s^F_{-2} \sim \cL_2 + (\cL_1)^2$, and in general
$$
\s_{-k}^F \sim
\sum (\cL_k)^{\a_k} (\cL_{k-1})^{\a_{k-1}} \cdots (\cL_1)^{\a_1}
$$
with the summation taken over $0 \leq \a_i \leq k$, for every $i$ and
$\a_1 + 2 \a_2 + \cdots + k \a_k = k$.  The proof is complete.
\end{proof}

%%%%%%%%%%%%%%%%%%%%%%%%%%%%%%%%%%%%%%%%%%%%%%%%%%%%%%%%%%%%%%%%%%

%%%%%%%%%%%%%  REFERENCES   %%%%%%%%%%%%%%%%%%%

%%%%%%%%%%%%%%%%%%%%%information about the author
\small
\vskip 1pc
{\obeylines
\noindent William \ J. \ Ugalde
\noindent Department of Mathematics
\noindent 14 MacLean Hall
\noindent The University of Iowa
\noindent Iowa City, Iowa, 52242
\noindent E-mail: wugalde@math.uiowa.edu
}

\end{document}